\magnification1200
\input amstex
\documentstyle{amsppt}
\def\M{\Cal M}
\def\K{\Cal K}
\def\w{\omega}
\def\ti{\times}
\def\a{\alpha}
\def\dlim{\varinjlim}
\def\IN{\Bbb N}
\def\IR{\Bbb R}
\def\s{\sigma}
\NoBlackBoxes

\vsize=19cm

\topmatter
\title
Topological classification of zero-dimensional $\Cal M_\w$-groups
\endtitle
\author
Taras Banakh
\endauthor
\address
Department of Mathematics, Lviv National University, Universytetska 1,
79000, Lviv, Ukraine
\endaddress
\email tbanakh\@franko.lviv.ua\endemail
\abstract
 A topological
group $G$ is called an $\M_\w$-group if it admits a countable cover
$\K$ by closed metrizable subspaces of $G$ such that a subset $U$ of $G$
is open in $G$ if and only if $U\cap K$ is open in $K$ for every $K\in\K$.

It is shown that any two non-metrizable
uncountable separable zero-dimenisional $\M_\w$-groups are homeomorphic.
Together with Zelenyuk's classification of countable $k_\w$-groups this
implies that the topology of a non-metrizable zero-dimensional
$\M_\w$-group $G$ is
completely determined by its density and the compact scatteredness
rank $r(G)$ which, by definition, is equal to the least upper bound of
scatteredness indices of scattered compact subspaces of $G$.
\endabstract
\thanks Research supported in part by grant INTAS-96-0753.
\endthanks
\endtopmatter

\document
In \cite{Ze} (see also \cite{PZ, \S4.3})
E.Zelenyuk has proven that the topology of a countable
topological $k_\w$-group $G$ is completely determined by its compact
scatteredness rank $r(G)$ which, by definition, is equal to the least upper
bound of scatteredness indices of compact scattered subsets of $G$. In this
note we extend this Zelenyuk's classification result onto the class of
punctiform $\M_\w$-groups.

Let us recall that a topological space $X$ is {\it scattered} if every
non-empty subset of $X$ has an isolated point. For a scattered space $X$ its
scatteredness index $i(X)$ is defined as the smallest ordinal $\a$ such that
the $\a$-th derived set $X^{(\a)}$ of $X$ is finite. Derived sets
$X^{(\beta)}$ of $X$ are defined by transfinite induction: $X^{(0)}=X$,
$X^{(1)}$ is the set of all non-isolated points of $X$;
$X^{(\beta+1)}=(X^{(\beta)})^{(1)}$ and
$X^{(\beta)}=\bigcap_{\gamma<\beta}X^{(\gamma)}$ if $\beta $ is a limit
ordinal. It can be easily shown that $i(X)<\w_1$ if $X$ is a hereditarily
Lindel\"of scattered topological space (in particular, a countable
compactum). For a topological space $X$ let
$$
r(X)=\sup\{i(K):K\text{ is a compact scattered subset of }X\}
$$
be the {\it compact scattered rank} of $X$.

A topological space $X$ is defined to be a {\it $k_\w$-space} (resp. an {\it
$\M_\w$-space}) if $X$ admits a countable cover $\K$ by compact Hausdorff
subspaces (resp. by closed metrizable subspaces) of $X$ such that a subset
$U$ of $X$ is open in $X$ if and only if $U\cap K$ is open in $K$ for every
$K\in\K$. A space $X$ is called an {\it $\M\K_\w$-space} if $X$ is both
a $k_\w$-space and an $\M_\w$-space.
A topological group $G$ is called a {\it $k_\w$-group} (resp.
{\it $\M\K_\w$-group,  $\M_\w$-group}) if its underlying topological space is
$k_\w$-space (resp. an $\M\K_\w$-space, an $\M_\w$-space). Since each
countable compactum is metrizable, we
conclude that each countable $k_\w$-space is an $\M\K_\w$-space.
On the other hand, according to Theorem 4 of \cite{Ba}, every
non-metrizable $\M_\w$-group is homeomorphic to the product $H\ti D$,
where $H$ is an open $\M\K_\w$-subgroup in $G$ and $D$ is a discrete space.

Following \cite{En$_2$, 1.4.3}, we say that a topological space $X$ is {\it
punctiform} if it contains no connected compact subspace containing more than
one point. Each punctiform $\sigma$-compact space is zero-dimensional
\cite{En$_2$, \S1.4}. On the other hand, there exist strongly infinite-dimensional
separable complete-metrizable punctiform spaces \cite{En$_2$,
6.2.4}.
Given a topological space $X$ by $d(X)$ its density is denoted.

\proclaim{Main Theorem} The topology of a non-metrizable punctiform
$\M_\w$-group is completely determined by its density and its compact
scatteredness rank. In other words, two non-metrizable punctiform
$\M_\w$-groups $G$, $H$ are homeomorphic if and only if $d(G)=d(H)$ and
$r(G)=r(H)$.
\endproclaim

To prove this theorem we need to make first some preliminary work. We say
that a topological space $X$ carries the direct limit topology with respect
to a tower $X_1\subset X_2\subset X_3\subset\dots$ of subsets of $X$ (this is
denoted by $X=\dlim X_n$) if $X=\bigcup_{n=1}^\infty X_n$ and a subset
$U\subset X$ is open if and only if $U\cap X_n$ is open in $X_n$ for every
$n\in\IN$.

Since the union of any two compact (resp. closed metrizable) subspaces in a
topological space is compact (resp. closed and metrizable, see \cite{En$_1$,
4.4.19}), we get the following

\proclaim{Lemma 1} A topological space $X$ is an $\M_\w$-space (an
$\M\K_\w$-space) if and only if $X$ carries the direct limit topology with
respect to a tower $X_1\subset X_2\subset\dots$ of closed metrizable
(compact) subsets of $X$.
\endproclaim

Under a {\it Cantor set} we understand a zero-dimensional metrizable
compactum without isolated points.

\proclaim{Lemma 2 \cite{Ke, 6.5}} Each uncountable metrizable compactum
contains a Cantor set.
\endproclaim

According to a classical theorem of Brouwer  \cite{Ke, 7.4}, each Cantor set is
homeomorphic to the Cantor cube $2^\w=\{0,1\}^\w$. It is well known that the
Cantor cube is universal for the class of metrizable zero-dimensional
compacta. In fact, it is universal is a stronger sense, see \cite{vE},
\cite{Po}.

\proclaim{Lemma 3} Suppose $A$ is a closed subset of a zero-dimensional
metrizable compactum $B$. Every embedding $f:A\to 2^\w$ such that $f(A)$ is
nowhere dense in $2^\w$ extends to an embedding $\bar f:B\to 2^\w$.
\endproclaim

Given a cardinal $\tau$ denote by $(2^\tau)^\infty=\dlim (2^\tau)^n$
the direct limit of the tower
$$
2^\tau\subset(2^\tau)^2\subset(2^\tau)^3\subset\dots
$$
consisting of finite powers of the Cantor discontinuum $2^\tau$
(here $(2^\tau)^n$ is identified
with the subspace $(2^\tau)^n\ti\{*\}$ of $(2^\tau)^{n+1}$, where $*$ is any
fixed point of $2^\tau$).

Using Lemma 3 by standard ``back-and-forth'' arguments (see \cite{Sa}) one
may prove

\proclaim{Lemma 4} A space $X$ is homeomorphic to $(2^\w)^\infty$ if and only
if $X$ is a zero-dimensional $\M\K_\w$-space satisfying the following property:
{\parindent30pt
\item{($\Cal S\Cal U$)} every embedding $f:B\to X$ of a closed subspace $B$
of a zero-dimensional metrizable compactum $A$ may be extended to an
embedding $\bar f:A\to X$.
\item{}}
\endproclaim

Now we are able to prove a ``separable'' version of Main Theorem.

\proclaim{Theorem} Every non-metrizable uncountable separable punctiform
$\M_\w$-group is homeomorphic to $(2^\w)^\infty$.
\endproclaim

\demo{Proof} Suppose $G$ is a non-metrizable uncountable separable punctifurm
$\M_\w$-group. It follows from Theorem 4 of \cite{Ba} that $G$ is an
$\M\K_\w$-group. Then $G$, being $\s$-compact and punctiform, is
zero-dimensional, see \cite{En$_2$, \S1.4}. According to Lemma 4, to show
that $G$ is homeomorphic to $(2^\w)^\infty$ it remains to verify the
property $(\Cal S\Cal U)$ for the group $G$.

Fix any embedding $f:B\to G$ of a closed subspace of a metrizable
zero-dimensional compactum $A$. By the continuity of the multiplication $*$
on $G$, the set $f(B)^{-1}*f(B)=\{f(b)^{-1}*f(b'):b,b'\in B\}\subset G$
is compact. It
follows from Theorem 4 of \cite{Ba} that there exists a sequence
$(x_n)_{n=1}^\infty\subset G$ converging to the neutral element $e$ of $G$
and such that $x_n\notin f(B)^{-1}*f(B)$ for every $n\in\IN$. This implies
that $f(B)$ is a nowhere dense subset in the compactum $f(B)*S_0$, where
$S_0=\{e\}\cup\{x_n:n\in\IN\}$. Next, since the $\M\K_\w$-group $G$ is
uncountable and $\s$-compact, it contains an uncountable metrizable compactum
which in its turn, contains a Cantor set $C\subset G$ according to Lemma 2.
Without loss of generality, $C\ni e$. It can be easily shown that the
compactum $f(B)*S_0*C$ has no isolated point and contains $f(B)$ as a nowhere
dense subset. Since $f(B)*S_0*C$ is a zero-dimensional metrizable compactum
without isolated points, it is homeomorphic to the Cantor cube $2^\w$, which
allows us to apply Lemma 3 to produce an embedding $\bar f:A\to
f(B)*S_0*C\subset G$ extending the embedding $f$. Thus the space $G$
satisfies the condition $(\Cal S\Cal U)$ and $G$ is homeomorphic to
$(2^\w)^\infty$.\qed
\enddemo

\proclaim{Lemma 5} If $G$ is a non-metrizable $\M_\w$-group, then
$r(G)\le\w_1$. Moreover, $r(G)=\w_1$ if and only if $G$ contains a Cantor set.
\endproclaim

\demo{Proof} Suppose $G$ is a non-metrizable $\M_\w$-group. Write $G=\dlim
M_i$, where $M_1\subset M_2\subset \dots$ of a tower of closed metrizable
subspaces of $G$ with $G=\bigcup_{i=1}^\infty M_i$. It follows that each
scattered compactum $K\subset G$ is contained in some $M_i$ and being
metrizable and scattered, is countable, see Lemma 2. Consequently,
$r(K)<\w_1$ for every such $K\subset G$. Hence $r(G)\le\w_1$.

If $G$ contains a Cantor set $C$, then $r(G)\ge r(C)\ge\w_1$ because $C$,
being universal in the class of zero-dimensional metrizable compacta,
contains copies of all countable compacta (whose scatteredness indices run
over all countable ordinals, see \cite{Ke, 6.13}).

Assume finally that $r(G)=\w_1$. According to Theorem 4 of \cite{Ba}, $G$ is
homeomorphic to the product $H\ti D$ of an $\K\M_\w$-group $H\subset G$ and a
discrete space $D$. Clearly, $\w_1=r(G)=r(H\ti D)=r(H)$. Write $H=\dlim K_i$,
where $K_1\subset K_2\subset\dots$ is a tower of metrizable compacta in $H$.
One of these compacta is uncountable (otherwise we would get
$r(H)=\sup\{r(K_i):i\in\IN\}<\w_1$, a contradiction with $r(H)=\w_1$).
Consequently, the group $H$ contains a Cantor set $C$, see Lemma 2.\qed
\enddemo

\demo{Proof of Main Theorem} Suppose $G_1$, $G_2$ are two non-metrizable
$\M_\w$-groups with $r(G_1)=r(G_2)$ and $d(G_1)=d(G_2)$. By Theorem 4 of
\cite{Ba}, for every $i=1,2$ the space $G_i$ is homeomorphic to the product
$H_i\ti D_i$, where $H_i\subset G_i$ is an $\K\M_\w$-group and $D_i$ is a
discrete space. Since $d(G_1)=d(G_2)$ and the spaces $H_1, H_2$ are separable,
we may assume that $|D_1|=|D_2|$ (if $d(G_1)=d(G_2)$ is countable, then
replacing $H_i$ by $G_i$, we may assume that $|D_1|=|D_2|=1$). Thus to prove
that the groups $G_1$ and $G_2$ are homeomorphic, it suffices to verify that
the groups $H_1$ and $H_2$ are homeomorphic. Observe that $r(G_i)=r(H_i\ti
D_i)=r(H_i)$ for $i=1,2$ and hence $r(H_1)=r(H_2)$.

If $r(H_1)=r(H_2)<\w_1$, then by Lemmas 2 and 6, the $\K\M_\w$-groups $H_1$ and
$H_2$ are countable and by Zelenyuk's theorem \cite{Ze}, they are
homeomorphic. If $r(H_1)=r(H_2)=\w_1$, then we may apply Theorem and
 Lemmas 2, 5 to
conclude that both groups $H_1$ and $H_2$ are homeomorphic to $(2^\w)^\infty$.
\qed
\enddemo

A topological space $X$ is defened to be an {\it
AE(0)-space} if every continuous map $f:B\to X$ from a closed subset of a
zero-dimensional compact Hausdorff space $A$ can be extended to a
continuous map $\bar f:A\to X$.

\proclaim{Conjecture} An uncountable zero-dimensional $k_\w$-group $G$
is homeomorphic to $(2^\tau)^\infty\times 2^\kappa$ for some cardinals
$\tau\le\kappa$ if and only if $G$ is an AE(0)-space.
\endproclaim

\enddocument